\documentclass[a4paper,12pt,reqno]{amsart}
\usepackage{amssymb}\theoremstyle{plain}
\usepackage{bm}

\usepackage{amscd}
\usepackage{enumitem}

\newtheorem*{theorem*}{Theorem}

\newtheorem*{corollary*}{Corollary}

\newtheorem*{lemma*}{Lemma}
\newtheorem{proposition}{Proposition}
\newtheorem*{proposition*}{Proposition}

\newtheorem*{conjecture*}{Conjecture}
\theoremstyle{definition}
\newtheorem{definition}{Definition}
\newtheorem*{definition*}{Definition}
\theoremstyle{remark}
\newtheorem{remark}{Remark}
\newtheorem*{remark*}{Remark}
\newtheorem{example}{Example}
\newtheorem*{problem*}{Problem}

\begin{document}

\begin{center}
\title[QM RZ LS AE]{Quantum models of the Riemann zeta function, lattice spin models and algebraic models of entanglement}
\end{center}

\maketitle

\begin{center}
{\bf Nikolaj M. Glazunov } \end{center}

\begin{center}
{\rm Glushkov Institute of Cybernetics NASU, Kiev, } \\

%\orcid{Orcid number}
{\rm  Email:} {\it glanm@yahoo.com }
\end{center} 

\bigskip


\begin{abstract}
A brief overview of results concerning the connection between the Hilbert-Polya  conjecture and the Riemann hypothesis about the Riemann zeta function, some new results on p-adic quantum computing, quantum entanglement based on lattice spin models and algebraic entanglement models is given. Quantum computing uses both photons and electrons, so their known properties are (very briefly) presented.

%The content of the paper is presented in more detail in the table of Contents and in the Introduction.
\end{abstract}

\tableofcontents

\section{Introduction}
The Hilbert-Pólya conjecture relates the Riemann Hypothesis about the Riemann zeta
 function to the eigenvalues of quantum  Hamiltonians. 
This also applies to Selberg`s trace formula and Connes approach based on noncommutative geometry. The corresponding problems and related selected current problems in quantum mathematics will be briefly reviewed. The examination is planned to be conducted on the basis of spectral methods and the adelic approach.
   Continuing my talk at the Seminar on December 30, 2025, an overview of new (2026) results in $p$-adic quantum computing will be presented.
   Quantum entanglement plays an important role in quantum information theory and quantum computing. Quantum entanglement based on lattice spin models will be presented.
   Due to the problems in understanding the physical basis of quantum entanglement, researchers introduce and investigate algebraic models of entanglement. A class of algebraic models of entanglement based on the theory of elliptic curves and Galois representations is presented.

\section{Quantum Hamiltonians}  
  \subsection{International System of Units, SI}

  metre m ,    kilogram kg,  second s  (sec)\\

 newton (symbol: N)   $N = \frac{kg\cdot m}{s^{-2}}$ the force that accelerates a mass of one kilogram at one metre per second squared 

  joule (symbol: J) (Newton meter)     is the unit of energy in SI:    
$  J =  \frac{kg\cdot m^2}{s^{-2}}$

 $erg = 10^{-7} J $\\

 Planck's constant $h = 6.6 \cdot 10^{-27} erg\cdot sec =
 6.6 \cdot 10^{-27}\cdot 10^{-7} J\cdot sec  = 6.6 \cdot 10^{-34 } J \cdot sec \;
(Newton \; meter\cdot sec)$\\

Quantum of action (Planck):    $h'= \frac{h}{2\pi}$            \\

Planck length: $\ell_P  \approx 1.6 \times 10^{-35} m$

 \subsection{Quantum Hamiltonian}
\subsubsection{de Broglie wave function}
\smallskip
 Photon:  speed of light $c$, oscillation frequency $\nu$,  wavelength
 $\lambda = \frac{c}{\nu}$, $\omega = 2\pi \nu$,  wave vector ${\mathbf k}, 
|{\mathbf k}| = \frac{2\pi}{\lambda}$ \\
 {}
Momentum and energy of a photon: ${\mathbf p} = h' {\mathbf k}$, $E = h' \omega$.\\

{\it Wave function}
\begin{equation}
\label{wf}
\psi = e^{\frac{i}{h'} (x p_x + y p_x + z p_z - Et)} 
\end{equation}

 {\it Pure quantum states:}\\
A system is in a pure state when it is possible to present it as a single (ket) vector 
 in a Hilbert space: 
 \(\vert \psi \rangle\) .

\subsubsection{Hamiltonians}
 {\it Hamiltonian for a photon:} 
\begin{equation}
\label{hp}
{\hat H} = h' \omega (a^{\dagger} a + \frac{1}{2})
\end{equation}
$a^{\dagger}_n$ adds a photon to the mode (photon criation operator)
\begin{equation}
\label{pco}
 (a^{\dagger}_n f)(x) = x_n f(x))
\end{equation}
$a_n$ removes a photon (photon annihilation operator)\\
 \begin{equation}
\label{pao}
(a_n f)(x) = \frac{\partial f}{\partial x_n}(x)
\end{equation}
 {\it Operators for the components of momentum} \\
\begin{equation}
\label{om}
 \hat P_x = ih'\frac{\partial}{\partial x}, \; \hat P_y = ih'\frac{\partial}{\partial y}, \; 
\hat P_z = ih'\frac{\partial}{\partial z},
\end{equation}
 The action of operators (\ref{om}) on the wave function (\ref{wf}) yields momentum  components that are the eigenvalues of these operators:
\begin{equation}
\label{ow[}
\hat P_x \psi = p_x \psi, \, p_x, \, \hat P_y \psi = p_y \psi, \, 
\hat P_z \psi = p_z \psi, 
\end{equation}  

 {\it Hamiltonian} 
%(according to Schr\"{o}dinger). 
The energy operator of a particle of mass $m$ in a field with potential energy $U$:
\begin{equation}
\label{ho}
\hat H = \frac{1}{2m}(\hat P_x^2  +  \hat P_y^2 + \hat P_z^2) + U(x,y,z)
\end{equation}
or
\begin{equation}
\label{hl}
\hat H = -\frac{{h'}^2}{2m} \left(\frac{{\partial}^2}{\partial x^2} + \frac{{\partial}^2}{\partial y^2}  + \frac{{\partial}^2}{\partial z^2} \right) + U(x,y,z) = 
 -\frac{{h'}^2}{2m} \triangle + U
\end{equation}
 where   $\triangle$ is the Laplacian.

Energy operator and Schr\"{o}dinger equation
\begin{equation}
\label{sche}
  -\frac{{h'}^2}{2m} \triangle \psi + U \psi = E \psi
\end{equation}
 
\subsection{Operators of orbital angular  momentum of an electron}
\label{ooan}
\begin{remark}
Operators of proper angular momentum of an electron are given in subsection \ref{lsm}.
\end{remark}
 \subsubsection{Orbital angular momentum of an electron}
  \begin{equation}
 \label{oam}
    \hat M_x = y \hat P_x - x \hat P_y, \:  \hat M_y = z \hat P_x - x \hat P_z, \:
  \hat M_z = x \hat P_y - y \hat P_z, 
 \end{equation}

These operators satisfy the following commutation relations:
 \begin{equation}
 \label{ocr}
\begin{array}{rcl}
 \hat M_y  \hat M_z -  \hat M_z  \hat M_y  = i h' \hat M_x, \\
 \hat M_z  \hat M_x -  \hat M_x  \hat M_z  = i h' \hat M_y,  \\
  \hat M_x  \hat M_y -  \hat M_y  \hat M_x  = i h' \hat M_z.
\end{array}
\end{equation}
 \begin{proposition}
  Let 
${\hat {\mathbf M}}^{(1)} = ({\hat M_x}^{(1)}, {\hat M_y}^{(1)},  {\hat M_z}^{(1)})$
 and
 ${\hat {\mathbf M}}^{(2)} = ({\hat M_x}^{(2)}, {\hat M_y}^{(2)},  {\hat M_z}^{(2)})$
be commuting operators 
${\hat {\mathbf M}}^{(1)} {\hat {\mathbf M}}^{(2)}  = 
{\hat {\mathbf M}}^{(2)} {\hat {\mathbf M}}^{(1)}$
that satisfy relations (\ref{ocr}). Then their sum
 ${\hat {\mathbf M}}^{(1)} + {\hat {\mathbf M}}^{(2)}$ also satisfies commutation relations (\ref{ocr}).
\end{proposition}

 \section{Riemann Hypothesis}
\subsection{Riemann  zeta function and the Mellin transform of theta
 function} \cite{tit,kavo,kob,naka}
\subsubsection{Representation}
\begin{equation}
\label{rz}
  \zeta(s) = \sum_{n = 1}^{\infty} \frac{1}{n^s} = \prod_p (1 - \frac{1}{p^s})^{-1}
\end{equation}
  $s = \sigma + it$

   Let us remind you that:: 

Critical line  $\Re  s = \frac{1}{2}$.

Critical strip $0 \le \sigma \le 1$ in the complex plane.

Let $c > 0$ be any constant.
Let us recall that the gamma function $\Gamma$ is a multiple $c^s$ Mellin transform
${\mathcal M}(\cdot) 
= \int_{0}^{\infty} () t^s \frac{dt}{t}$
 of the function $e^{-ct}$:\\
  $ \int_{0}^{\infty} e^{-ct} t^s \frac{dt}{t}$.\\
   Let's perform a change of variables $t = \frac{u}{c}$.\\
 $c^s (\int_{0}^{\infty} e^{-ct} t^s \frac{dt}{t}) =
  \int_{0}^{\infty} e^{-u} u^s \frac{du}{u} =  \int_{0}^{\infty} e^{-u} u^{s-1} du = 
\Gamma(s)$.\\

  $\Gamma(\frac{s}{2}) = \int_0^{\infty} e^{-t} t^{\frac{1}{2}s - 1} dt$

\subsubsection{Functional equation}
   Let $t > 0$. Define the theta function
\begin{equation}
\label{theta}
  \theta(t) = \sum_{-\infty}^{\infty} e^{-\pi t n^2}.
\end{equation}
\begin{equation}
\label{xi}
  \xi(s) =  \frac{1}{2} s(s-1){\pi}^{-\frac{s}{2}} \Gamma(\frac{s}{2})\zeta(s) 
\end{equation}

\begin{equation}
\label{fe}
   \xi(s) = \xi(1 - s) 
\end{equation}

\begin{remark}
  Zeros of  the gamma function  $\Gamma(\frac{s}{2})$ is called  the trivial  zeros of 
 $\zeta(s)$.
\end{remark}

\begin{remark}
Equations (\ref{xi}) and (\ref{fe}) can be obtained by applying the Mellin transform to the theta function   (\ref{theta}).
\end{remark}

\subsubsection{The Riemann hypothesis}
The Riemann hypothesis (RH) states that the non-trivial complex zeros of  $\zeta(s)$ lie on the critical line $\Re  s = \frac{1}{2}$.

\subsubsection{On A. Speiser result} According to the results of A. Speiser, the Riemann's zeta function $\zeta(s)$ satisfies RH if and only if  $\zeta^{\prime}(s)$ has no zeros left of the critical line $\Re  s = \frac{1}{2}$.
\subsection{Quantum Hamiltonians to the Riemann hypothesis} 
 Operators whose eigenvalues correspond to the locations of the nontrivial zeros of the zeta function
\subsubsection{Berry-Keating, derived and  others Hamiltonians}
\begin{equation}
\label{bkh}
\hat H_{BK} =  \hat  x \hat p +    \hat p  \hat x
\end{equation}
\begin{equation}
\label{bybkh}
  \frac{1}{2}(1 - i \hat H_{BK} )
\end{equation}
\begin{equation}
\label{zbkh}
\zeta( \frac{1}{2}(1 - i \hat H_{BK} ))
\end{equation}

$(-i(\frac{1}{2}+x\frac{d}{dx}))$

\section{On Selberg approach and Selberg classes}
 \subsection{On Selberg`s trace formula}
  Let now $-\triangle$  be the Laplace-Beltrami operator (Quantum Hamiltonian)  defined on a matifold $\mathbf M$.
\begin{example}
   Simple one-dimensional case.\\
   $-\triangle = -\frac{d^2}{dx^2}$ on the circle ${\mathbb S}^1$ of length $2\pi$,\\
  The eigenvalues of $-\frac{d^2}{dx^2}$ are $m^2, \; m \in {\mathbb Z}$, with corresponding eigenfunctions $\varphi_m(x) = (2\pi)^{-\frac{1}{2}}e^{imx}$. 
In the case the Selberg`s trace formula  is the Poisson summation formula: \\
   Let $L = -\frac{d^2}{dx^2}$ be the linear operator acting on the $2\pi-$periodic function by \\
$[Lf](x) = \int_0^{2\pi} k(x,y)f(y)dy$.\\
 with kernel $k(x,y) =
 \sum_{m \in {\mathbb Z}} h(m) {\varphi_m(x)} {\overline {\varphi}_m(y)}$ \\
We have: $L\varphi_m = h(m)\varphi_m$\\
 The corresponding Poisson summation formula is: \\
 \begin{equation}
\label{stf}
Tr \: L = \sum_{m \in {\mathbb Z}} h(m) =
 \sum_{m \in {\mathbb Z}} \int_{n \in {\mathbb Z}} h(\rho)e^{2\pi in{\rho}} d{\rho}.
\end{equation}
\end{example}

\section{On Selberg classes}
\subsection{Dirichlet series and Selberg classes}
  Let $a_n$ be complex numbers, $s = \sigma + it$ be a complex variable, and $f(s) =
  \sum_{n = 1}^{\infty} \frac{a_n}{n^s}$ be a Dirichlet series.
\begin{definition} (by A. Selberg)
 The Selberg class $S$ consists of functions $f(s)$ satisfying the following axioms: \\
(1) Dirichlet series $f(s)$  absolutely convergent for $\sigma > 1$.\\
(2) (Analytic continuation) There exists an integer $m$ such that \\
$(s - 1)^mf(s)$ is an entire function of finite order.\\
(3) Functional equation (Extension of (\ref{fe})).\\
(4) (Ramanujan conjecture) For every $\epsilon > 0$, $a(n) = O(n^{\epsilon})$.\\
(5) Euler product. (Extension of (\ref{rz})).
\end{definition}
 Let $S^{\#}$ be the extended Selberg class   for the Selberg class $S$.
 The extended Selberg class  $S{\#}$ for the Selberg class $S$ by A. Selberg [Collected Papers, Vol. II, pp. 47-63. Springer-Verlag, Berlin (1991; Zbl ] is defined by abandoning the Ramanujan conjecture and the Euler product. It was introduced by J. Kaczorowski and  A. Perelli in \cite{kape}.

As Dirichlet series $F(s) \in S^{\#}$ satisfies a functional equation
\[
  \Phi(s) = \omega \overline{\Phi(1 - s)}  
\]
  where 
 $\Phi(s) = F(s)Q^s \prod_{j=1}^r\Gamma(\lambda_js + \mu_j)$,
 the invariant $d_F = 2 \sum_{j=1}^r \lambda_j$ is called the degree of $F(s)$.

 %   $\zeta(X,s) = \prod_{x \in X_{0}} \frac{1}{1 - N(x)^{-s}},$ 
Theorem 1. \cite{zdzfesc2025}. Let $F(s) \in S_0^{\#}$, $q_F \ge 2$, and 
$a_1 \ne 0$. Let $M(T)$ be the number of zeros of $F(s)$
 in the region $0 < t < T$ (counting with multiplicities). For every $\epsilon > 0$, there exists a sequence $\{T_j\}$ 
(dependent on $F(s)$), $T_j \to \infty$, such that
\[
  \vert M(T_j) -  \frac{\log q_F}{2\pi} T_j \vert \le \epsilon,
 \]
for all $j \in {\mathbb N}.$

Theorem 2.  \cite{zdzfesc2025}.   Let $N(T)$ and $N_1(T)$ denote the number of zeros (counting with multiplicities) in the region $\sigma < 1/2$ and $0 <t< T$ of $F(s) \in S^{\#}_0$, $q_F \ge 2$, and  $F^{\prime}(s)$ respectively. There exists a sequence $\{T_j\}$  (dependent on $F(s)$), $T_j \to \infty$, such that $N(T_j) = N_1(T_j)$, for all $j \in {\mathbb N}$.

Corollary 1.3. $F(s) \in S^{\#}_0$ satisfies RH if and only if $F^{\prime}(s)$ has no zeros in $\sigma < 1/2$.

Corollary 1.4. Let $F(s) \in S^{\#}_0$, $q_F \ge 2$. Let $A_n = a_n/\sqrt n$. If
$|A_{q_F}| \ge \frac{1}{2} \sum_{ n | q_F, \; n \ne 1,q_F} |A_n| $
is satisfied then $F(s)$ satisfies RH.

\section{On noncommutative spaces and Alane Connes (with his collaborators) works}
 "Zeta zeros and prolate wave operators : semilocal adelic operators"
The essence of this interesting paper is the integration in the semilocal framework by 
 Connes,  Consani \cite{coco}  
 of the trace formula by Connes [Selecta Math. (N.S.
 5, 1, 29-106, (1999; Zbl 0945.11015)] two recent results by the first and second authors
[ Enseign. Math., 69, 1-2, 93-148 (2023; Zbl 1537.11121)] and by the first and therd authors  [Proc. Natl. Acad. Sci. USA, 119, 22 (2022; Zbl 07998560) ]
``on the spectral realization of the zeros of the Riemann zeta function by introducing
a semilocal analogue of the prolate wave operator''. 

Based on the (quantum) mathematical physics by Connes, Marcolli \cite{coma}, the results of the two above-mentioned papers
are interpreted as the infrared (low-lying) part and the ultraviolet behavior of the zeros
 of the Riemann zeta function. 
 
 The main results of the paper  are Theorem 1 on semilocal Hardy–Titchmarsh transform   and Theorem 2  that relate with ``the sought-for Weil cohomology''.
  The paper  consists of six sections.
   The first section is of introductory nature and designed to acquaint the reader with main ideas and results of the paper.
Section 2 presents cyclic pairs and associated prolate operators.

Let $\mathcal H$ be a Hilbert space, $D$ a selfadjoint  operator acting on 
$\mathcal H$ and $\xi \in \cap_{n \in {\mathbf N}} Dom D^n$ a unit cyclic vector. Let $\lambda > 0$.
For a cyclic pair $(D, \xi)$ authors of the paper under review define the formal prolate operator $\omega(D,\xi,\lambda)$ as 
$\omega(D,\xi,\lambda) = -D^2 + \lambda^2 N$, where $N$  is the  number operator ($N{\xi}_j = j{\xi}_j$).

Let the ${\mathcal C}^*$-algebra 
${\mathcal A} := c_0(\mathbb N)$ be the algebra of sequences 
vanishing at $\infty$.
Authors for an even cyclic pair $(D, \xi)$ and the ${\mathbb Z}/2$ grading $\gamma$ define the even spectral triple $({\mathcal A}, {\mathcal H}, D)$   where the ${\mathcal C}^*$-algebra ${\mathcal A}$ acts in the Hilbert space ${\mathcal H}$.
After a number of preliminary results (Propositions 2.1 - 2.3),
 the authors define the Jacobi matrices of the formal prolate
operator (Proposition 2.4).

Section 3 deals with Hardy–Titchmarsh transform: archimedean place.
 ``Section 4 is devoted to the extension of the Hardy–Titchmarsh transform in the
semilocal situation involving a finite set S of places of Q containing the archimedean
one, and the analysis of the semilocal Sonin space.''
Section 5  studies metaplectic representation in the Jacobi picture.
  The sixth section is the Appendix: Fourier transform on $p$-adic numbers.
  The (Weil) positivity is given in the paper by A. Weil [Comm. Lund., pp. 252–267 (1952; Zbl 0049.03205)]. See also papers by 
 H. Yoshida [In: Zeta Functions in Geometry (Tokyo, 1990), Adv. Stud. Pure Math., vol. 21, pp. 281-325. Kinokuniya, Tokyo (1992; Zbl 0817.11041) ]  and by E Bombieri [Rend. Mat. Acc. Lincei, s. 9, v. 11:183-233 (2000; Zbl )].

\section{About spin models}

\subsection{About $p-$adic spin models}
 Recall at first about Planck length: $\ell_P  \approx 1.6 \times 10^{-35} m$

\begin{remark}  
\label{sps}
  "sub-Planckian" studies explore phenomena, structures, and theoretical limits well below or beyond these bounds.
\end{remark}
\begin{remark}  
\label{spsna}
The hypothesis that the geometry of spacetime is \\
{\bf non-Archimedean} at the sub-Planckian scales.
This motivates the study of non-Archimedean (in particular, $p-$adic) quantum mathematics.
\end{remark}

  In particular, to study qubits and spin in the non-Archimedean case, the $p-$adic rotation group $SO(3)_p$ and the components $SO(3)_p \; mod \: p^n$ of their projective limit are studied \cite{svampa}.

\subsection{Spin and quantum operators}
 % Let's continue our considerations  \ref{ooan}.

\section{About lattice spin models}

   Here we will very briefly touch on the Pólya enumeration theory in combinatorial analysis, its connection with the study of particle spin on lattices, and their generalizations \cite{ahst,mopo,pol}.

%   Consider 
  Let $\chi$ be a general spin state  for a single spin $1$ or $1/2$ particles.

$\Psi({\mathbf r}, spin) = \psi({\mathbf r})\chi$

$\chi = \left( \begin{array}{c} a \\
b \end{array} \right)$.

\section{On the random walk of a particle on a lattice}

Let $\Lambda $ be a  lattice with base 
\[
\{{\mathbf a}_1, \ldots ,{\mathbf a}_n \}
\] 
where all vectors  ${\mathbf a}_i$ have integral
coordinates \cite{Cassels}.

%Let us 
 At first consider the Pólya problem of a random walk of a particle on a lattice.

Let us consider a lattice $\Lambda_n $ of dimension $n$, the basis of which are the unit vectors ${\mathbf e}_1 = (1, 0 \ldots, 0), \ldots ,{\mathbf e}_n = (0, 0 \ldots, 1)$ of the real space ${\mathbf R}^n$.

We assume that there can be particles at the lattice points, and that particles can move between lattice points.
 \begin{example}
Let at the moment $t = 0$ the particle be at the zero point  (at the origin)  of a one-dimensional lattice  $\Lambda_1 = \{-\infty <m  < \infty, m \in {\mathbb Z} \}  \simeq {\mathbb Z}$ and can move one unit along the lattice to the left or right with a probability of $\frac{1}{2}$.
 \end{example}

 \begin{definition}
Let us call the probabilistic structure on the lattice for the possible movement of a particle by one step with probability $\frac{1}{2}$ in any of the admissible directions the expression \\
1) $\frac{1}{2}(e^{i \varphi} + e^{-i \varphi} )$ for a one-dimensional lattice,\\
2) $\frac{1}{4}(e^{i {\varphi}_1} + e^{-i {\varphi}_1} + e^{i {\varphi}_2} + e^{-i {\varphi}_2})$ for a two-dimensional lattice.
  \end{definition}

\begin{remark}
   In a similar way, one can define probabilistic structures for
lattices of dimension $n$ greater than 2.
 \end{remark}
\begin{remark}
Let us define for the lattice $\Lambda_1$ the integral operator $\frac{1}{2\pi} \int_{-\pi}^{\pi} (\cdot) d {\varphi}$ acting on the probabilistic structure of the one-dimensional lattice  (on corresponding exponential functions).
   As a result of this action, we obtain the coefficient of the term $e^{i \varphi m }$, which is the probability that after $m$ steps the particle will end up at point $m$.
 \end{remark}
 \begin{definition}
Let at the moment $t = 0$ the particle be at the zero point  (at the origin)  of the lattice  
 $\Lambda_n $.
Let $p_t({\mathbf m})$ be the probability that a particle will end up at point 
${\mathbf m}$ of the lattice after $t$ steps.
 \end{definition}
 
 \begin{example}
Let us consider the simplest case of a particle moving along a one-dimensional
 lattice $\Lambda_1.$\\
  For the one-dimensional lattice the probabilistic structure has the form \\
  $\frac{1}{2}(e^{i \varphi} + e^{-i \varphi} )$ ,\\
  $ 
p_1( 1) = \frac{1}{2\pi} \int_{-\pi}^{\pi} \left(\frac{1}{2}(e^{i \varphi} + e^{-i \varphi}) \right)^1 e^{-i \varphi} d {\varphi} = \frac{1}{2}.
 $
  \end{example}

\begin{proposition}
\label{pl}
  For an $n$-dimensional lattice $\Lambda_n $ with $2n$ possible equally probable movements from each point of the lattice to an adjacent point, we have \\
 \begin{equation}
\label{stf}
p_t({\mathbf m}) = \frac{1}{(2\pi)^n} \int_{-\pi}^{\pi} \cdots \int_{-\pi}^{\pi}
 \left(\frac{1}{2n} \sum_{j=1}^{n} (e^{i \varphi_j} + e^{-i \varphi_j} )\right)^t e^{-i {\mathbf m}\cdot \bm{ \varphi}} d {\varphi}_1 \ldots d {\varphi}_n.
\end{equation}
  = 
 \begin{equation}
\label{stfr}
\frac{1}{(2\pi)^n} \int_{-\pi}^{\pi} \cdots \int_{-\pi}^{\pi}
 \left(\frac{1}{n}(\cos \varphi_1 + \cdots   + \cos \varphi_n)\right)^t e^{-i {\mathbf m}\cdot \bm{ \varphi}} d^n  \bm{ \varphi}
\end{equation}
  where
$
     {\mathbf m}\cdot \bm{ \varphi} = m_1 {\varphi}_1 + \cdots + m_n {\varphi}_n
$
   and 
$d^n  \bm{ \varphi} = d {\varphi}_1 \ldots d {\varphi}_n$
\end{proposition}

\subsection{On lattice spin models}
\label{lsm}
\subsubsection{Spin (the electron's own angular momentum) and Pauli matrices. }
 Let's continue our considerations  \ref{ooan}.
   According to Pauli and Dirac \cite{dirac,fock}, the angular momentum 
${\hat {\mathbf M}}$ of an electron is the sum of the orbital momentum
 ${\hat {\mathbf M}}_{orb}$ and the proper momentum ${\hat {\mathbf M}}_{pro}$:
$$
  {\hat {\mathbf M}} = {\hat {\mathbf M}}_{orb} + {\hat {\mathbf M}}_{pro}.
$$

The square ${\hat {\mathbf M}}_{pro}^2$ of the electron's proper moment takes only one value, equal to ${h'}^2s(s + 1)$ at $s = \frac{1}{2}$ (electron spin).\\
 Pauli matrices 
 \begin{equation}
  \sigma_x =  \left( \begin{array}{cc} 0 & 1 \\ 
 1 & 0  \end{array} \right),
  \sigma_y = \left( \begin{array}{cc} 0 & -i \\ 
 i & 0  \end{array} \right),
\sigma_z = \left( \begin{array}{cc} 1 & 0 \\ 
 0 & -1  \end{array} \right)
 \end{equation}
 \begin{equation}
 \label{pro} 
 (\hat M_x)_{pro} =\frac{{h'}}{2} \sigma_x, \:  (\hat M_y)_{pro} =\frac{{h'}}{2} \sigma_y, \:
   (\hat M_z)_{pro} =\frac{{h'}}{2} \sigma_z, 
 \end{equation}\\

   The wave function of an electron can be represented as 
$\psi = \psi({\mathbf r}, \sigma)$, where ${\mathbf r}$ is the radius vector of the electron and $\sigma$ is the spin variable (spin direction), which takes two values: $1$ (spin up 
$\uparrow$) and $-1$ (spin down $\downarrow$).

\subsubsection{On Ising model. Spins plased on points of $\Lambda$.}
 In the Ising model, each point of the lattice $\Lambda$ is assigned a variable 
$\sigma_j$, which denotes the direction of spin at that point.
\subsubsection{On Ising model.  Spins   interact with their neighbors}
Let's introduce the interaction energy $E_s$ of adjacent spins. It depends on the orientation of these spins.\\
If the spins are oriented parallel, i.e., $\uparrow \uparrow$, then the interaction energy is $-E_s$. \\
If the spins are oriented antiparallel, i.e., $\uparrow \downarrow$, then the interaction energy is $E_s$.
\begin{remark}
   Interaction energy  of adjacent spins is equal $-E_s \sigma \sigma'$. \\
   The position of minimum energy of a pair of spins $\sigma$ and $\sigma'$ is the position in which they are parallel.
\end{remark}

%Geometric lattice $\Lambda$.

\section{On entanglement in algebraic extensions of algebraic number fields} 
\subsection{About quantum entanglement}
A quantum state is entangled  if its density matrix can not be written as a convex sum
 of product states.
 \subsection{On entanglement in algebraic extensions}
 \begin{definition*}\cite{campen}.
Let $K$ be a field and let ${\mathcal K} = \{K_n\}$ a family of Galois extensions of $K$ inside an algebraic closure of $K$. The family ${\mathcal K}$ is called  linearly disjoint over $K$ if for the compositum $F$ of the fields $K_n$, the natural inclusion map
\begin{equation}
   Gal(F/K) \hookrightarrow \prod_{n \in M} Gal(K_n/K)
\end{equation}
is an isomorphism. If this is not the case, then the family ${\mathcal K}$F is called entangled over $K$.
 \end{definition*}

This phenomenon was first discovered by Serre.

This situation (entanglement) appearing in the family of division fields of elliptic curves with complex multiplication. \\

{\bf Acknowledgments.} 
I would like to thank Academician V. Gusynin for his consultations.
I thank the leaders and the scientific secretary of the International scientific seminar "Quantum Computing" for the invitation to give a talk.
   I would like to thank Academician V. Zadiraka for his questions and useful discussions with him during my presentation.

%\huge{Thank you for your attention}
%{\bf ORCID}\\
%https://orcid.org/0000-0002-1586-2696 \\

\end{document}